\input amstex
\loadbold

\define\Pee{{\Bbb P}}

\define\Cee{{\Bbb C}}

\define\proof{\demo{Proof}}
\define\endproof{\qed\enddemo}

\define\theorem#1{\proclaim{Theorem #1}}
\define\lemma#1{\proclaim{Lemma #1}}
\define\proposition#1{\proclaim{Proposition #1}}
\define\corollary#1{\proclaim{Corollary #1}}
\define\claim#1{\proclaim{Claim #1}}

\define\section#1{\specialhead #1 \endspecialhead}
\define\ssection#1{\medskip\noindent{\bf #1}}

\documentstyle{amsppt}
\leftheadtext{}
\rightheadtext{}

\pageno=1
\topmatter
\title Composition law and Nodal genus-$2$ curves in $\Pee^2$
\endtitle
\author 
Sheldon Katz$^1$ $\quad$ Zhenbo Qin$^2$ $\quad$ Yongbin Ruan$^3$
\endauthor
\address 
Department of Mathematics, Oklahoma State University, 
Stillwater, OK 74078
\endaddress
\email   
katz\@math.okstate.edu
\endemail
\address 
Department of Mathematics, Oklahoma State University, 
Stillwater, OK 74078
\endaddress
\email   
zq\@math.okstate.edu
\endemail
\address Department of Mathematics, University of Wisconsin, 
Madison, WI 53706
\endaddress
\email ruan\@math.utah.edu 
\endemail
\thanks ${}^1$Partially supported by an NSF grant and an NSA grant
\endthanks
\thanks ${}^2$Partially supported by an NSF grant
\endthanks
\thanks ${}^3$Partially supported by an NSF grant and a Sloan fellowship
\endthanks
\endtopmatter

\TagsOnRight
\NoBlackBoxes
\document
\section{1. Introduction}

Enumerative algebraic geometry is an old field of algebraic geometry. 
There are many fascinating problems going back more than a hundred years
to the Italian school.  The most famous one is perhaps 
the counting problem for the number of holomorphic curves in $\Pee^2$. 
There are in fact two different counting problems. 
Let $\tilde N_{g,d}$ be the number of irreducible, reduced, nodal, degree-$d$ 
genus-$g$ curves which pass through $3d + (g-1)$ general points in $\Pee^2$. 
The integer $\tilde N_{g,d}$ is often referred to as the Severi number.
These numbers have recently been computed in \cite{C-H}.
A companion number is the number $N_{g,d}$ of irreducible, reduced, nodal, 
degree-$d$ genus-$g$ curves whose normalization has a fixed complex structure 
and which pass through $m(d)$ general points in $\Bbb P^2$. 
Here $m(d)$ stands for $3d -1$ when $g = 0, 1$ and $3d - 2(g-1)$ when $g \ge 2$.
Clearly, $\tilde{N}_{0,d}=N_{0,d}$. 

The goal of this paper is to compute $N_{2, d}$ for $d \ge 4$ (see (1.1)).
The lower degree cases for $N_{g,d}$ were known classically for many years. 
Not much progress has been made until the introduction of 
quantum cohomology theory. Based on ideas from physics, 
a recursion formula for $N_{0,d}$ was proved by \cite{R-T, K-M}. 
The influence of physical ideas opens up entirely new directions in
enumerative geometry. Roughly speaking, $N_{0,d}$ can be interpreted as 
a correlation function in a certain topological quantum field theory 
(topological sigma model \cite{Wit}).  All
topological quantum field theories have a composition law, which in this
instance gives the beautiful recursion formula for $N_{0,d}$. 
Furthermore, the composition law naively suggests a recursion formula 
for $N_{g,d}$ in terms of $N_{0,d}$. Unfortunately, a simple calculation of 
lower degree elliptic curves showed that the formula from physics 
always gives a wrong answer for higher genus case $g>0$. 
It was showed in \cite{R-T} that the correlation function $\Psi_{g,d}$ 
(Gromov-Witten invariants) of the topological sigma model counts 
the number of perturbed pseudo-holomorphic maps. Moreover, 
it satisfies the composition law physicists predicted and 
hence can be computed by $N_{0,d}$.
However, $\Psi_{g,d} \neq N_{g, d}$ for $g>0$. 
The original problem of computing $N_{g,d}$ remains to be solved. 
One obvious approach is to compute the error term $\Psi_{g,d}-N_{g,d}$ 
and then to use the formula of $\Psi_{g,d}$ to compute $N_{g,d}$. 
Such an approach involves some delicate obstruction analysis and 
was carried out in \cite{Ion} for $N_{1,d}$. At the same time, 
a direct argument for $N_{1,d}$ was given independently in \cite{Pa1}. 

To explain our approach, we have to explain the composition law of 
the topological sigma model. Recall that we fix a complex structure on 
a genus-$g$ curve $\Sigma_g$ to define $N_{g, d}$ (and $\Psi_{g,d}$). 
Roughly speaking, the composition law gives an explicit formula of 
$\Psi_{g,d}$ in terms of $\Psi_{g',d'}$ with $g'<g, d'\leq d$ 
when we degenerate $\Sigma_g$ to a stable curve. 
As we mentioned, such a composition law fails for $N_{g,d}$. 
Our observation is that if we degenerate $\Sigma_g$ to a stable curve $C_0$
with only rational components each of which contains exactly three nodal 
points, then an analogue of the composition law might still hold. 
In fact, we shall prove that for $d \ge 4$, 
$$N_{2, d} = {(d-1)(d-2)(d-3) \over 2d} N_d$$
$$+ \sum \limits_{d_1+d_2=d}
{d_1d_2(d_1d_2d-6d+18)-4d \over 12d}
{3d-2 \choose 3d_1 -1} d_1d_2 N_{d_1} N_{d_2} \eqno (1.1)$$
where for simplicity, we have used $N_d$ to stand for both $N_{0, d}$
and $\tilde N_{0, d}$. 
Our arguments are parallel to those of Pandharipande \cite{Pa1}.
We expect that the same method works for any $g$. 
The difficulty is a technical one which becomes harder as $g$ gets larger. 
However, we believe that $N_{g, d}$ is closely related
to the number of irreducible, reduced, degree-$d$ plane rational curves
which pass through certain general points in $\Pee^2$ and 
have certain types of singularities. Indeed, it is not difficult to see that
one term involved in $N_{g, d}$ with $g > 2$ is 
$$(3d - 2(g-1))! \cdot \sum \limits_{d_1 + \ldots + d_{2(g-1)} = d} \,
\prod \limits_{i=1}^{2(g-1)} {d_i^3 \cdot N_{d_i} \over (3d_i -1)!}.
\eqno (1.2)$$
It is interesting to note that this term appears explicitly in $\Psi_{g,d}$.
It would be very interesting to figure out the other terms in $N_{g, d}$.

This paper is organized as follows. In section 2, the formula (1.1) is proved. 
In the proof, we need to know the number of constraints on stable maps 
that are degenerations of maps on irreducible and smooth curves,
and the number of irreducible, reduced, degree-$d$ rational plane curves 
that pass through $(3d-2)$ general points in $\Pee^2$ and 
have exactly one triple point with all other singularities being nodes.
These two numbers are studied in section 3.

\section{2. Proof of (1.1)}

First of all, we recall some definitions and notations for $g \ge 2$. 
Let $\overline{\frak M}_g$ be the Deligne-Mumford moduli space of 
stable genus-$g$ curves, and let 
$$\overline{\frak M}_{g}(d) \quad {\overset \text{def} \to =} \quad 
\overline{\frak M}_{g, 3d-2(g-1)}(\Bbb P^2, d) \eqno (2.1)$$ 
be the moduli space of stable maps from $(3d - 2(g-1))$-pointed genus-$g$ 
curves to $\Bbb P^2$ such that the homology class of the images is 
$d[\ell]$ where $[\ell]$ stands for the homology class of a line $\ell$ 
in $\Bbb P^2$. Then, there is a natural map 
$\pi: \overline{\frak M}_{g}(d) \to \overline{\frak M}_g$
obtained by forgetting the stable maps and all the marked points, 
and then contracting any unstable components.
For a stable genus-$g$ curve $C$, let $\overline{\frak M}_C(d)$ be 
the moduli space of stable maps to $\Bbb P^2$ from curves $D$ 
stably equivalent to $C$ with $(3d - 2(g-1))$ marked points
such that the homology class of the images is $d[\ell]$. In the above, 
we say that $D$ is {\it stably equivalent\/} to $C$ if contracting 
the unstable components of $D$ yields $C$. By the universal properties of 
moduli spaces, there is a canonical bijection 
$$\overline{\frak M}_C(d) \to \pi^{-1}([C]) \eqno (2.2)$$
which is an isomorphism when $[C] \in \overline{\frak M}_g$ is general. 
Let $W_g(d) \subset \overline{\frak M}_{g}(d)$ be the locus of stable maps 
whose domains are irreducible,
and let $\overline{W}_g(d)$ be the closure of $W_g(d)$ in 
$\overline{\frak M}_{g}(d)$. Then $W_g(d)$ is a reduced and irreducible
open subset of dimension $6d - (g -1)$. 
For $i = 1, \ldots, 3d - 2(g-1)$, define the evaluation map
$e_i: \overline{W}_g(d) \to \Bbb P^2$ by 
$[\mu: (D, p_1, \ldots, p_{3d - 2(g-1)})] \mapsto \mu(p_i)$, and
$\Cal L_i = e_i^*(\Cal O_{\Bbb P^2}(1))$. Put 
$Z = c_1(\Cal L_1)^2 \cap \ldots \cap c_1(\Cal L_{3d-2(g-1)})^2\cap 
[\overline{W}_g(d)]$.  In the sequel, we will usually think of
$Z$ as the cycle determined by the condition that $\mu(p_i)$ is a fixed
general point of $\Pee^2$.
Now for a general curve $[C] \in \overline{\frak M}_g$, the intersection
$$\pi^{-1}([C]) \cap [\overline{W}_g(d)-{W}_g(d)]$$
has codimension at least one in the irreducible and reduced subvariety 
$\pi^{-1}([C]) \cap \overline{W}_g(d)$. 
Since the linear series $e_i^*|\ell|$ are base-point-free,
$$\pi^{-1}([C]) \cap Z = [\pi^{-1}([C]) \cap \overline{W}_g(d)] \cap Z
= [\pi^{-1}([C]) \cap {W}_g(d)] \cap Z.$$
Moreover, by Bertini's Theorem, the intersection cycle 
$[\pi^{-1}([C]) \cap {W}_g(d)] \cap Z$ consists of 
finitely many reduced points in $\pi^{-1}([C]) \cap {W}_g(d)$. 
The number of these points is precisely $N_{g, d}$.
Thus for a general curve $[C] \in \overline{\frak M}_g$,
$$N_{g, d} = [\pi^{-1}([C]) \cap {W}_g(d)] \cap Z
= \pi^{-1}([C]) \cap Z = \pi^{-1}([C]) \cdot Z.$$
It follows that for every stable curve $[C] \in \overline{\frak M}_g$, we have
$$N_{g, d} = \pi^{-1}([C]) \cdot Z. \eqno (2.3)$$

Next, we construct a special stable genus-$g$ curve $C_{0, g}$ 
in $\overline{\frak M}_g$ by induction on $g$. First of all, 
$C_{0, 2}$ consists of two smooth rational curves intersecting 
transversely at three points. To get $C_{0, 3}$, 
we blow-up two nodal points in $C_{0, 2}$ by adding two smooth rational curves 
which intersect transversely at one point. In general, 
to obtain $C_{0, g}$ from $C_{0, g-1}$, we blow-up two nodal points in 
$C_{0, g-1}$ by adding two smooth rational curves 
which intersect transversely at one point. 
Thus, $C_{0, g}$ consists of $2(g-1)$ smooth rational curves 
which are denoted by $R_1, \ldots, R_{2(g-1)}$, 
and has $3(g-1)$ nodal points. Moreover, if $g > 2$, 
then for each smooth rational curve $R_i$ in $C_{0, g}$, 
there exist three other smooth rational curves $R_{k_1}, R_{k_2}, R_{k_3}$ 
in $C_{0, g}$ such that $R_i$ and each $R_{k_j}$ ($j = 1, 2, 3$) 
intersect transversely at one point. 
For simplicity, we denote the curve $C_{0, g}$ by $C_0$. 
Let $\omega_{C_0}$ be the dualizing sheaf of $C_0$. 
Then the restriction of $\omega_{C_0}$ to each smooth rational curve $R_i$ 
in $C_0$ has degree $1$. Thus if $L$ is a line bundle on $C_0$ such that
$\text{deg}(L|_{R_i}) \ge 0$ for all $i$ with $1 \le i \le 2(g-1)$,
$\text{deg}(L|_{R_{i_1}}) > 1$ for at least one $i_1$, 
and $\text{deg}(L|_{R_{i_2}}) = 0$ for at most one $i_2$, then 
$$H^1(C_0, L) \cong H^0(C_0, L^{-1} \otimes \omega_{C_0}) = 0. \eqno (2.4)$$

Let $[C] \in \overline{\frak M}_g$ be generic,
and $|\text{Aut}(C)|$ be the order of the automorphism group of $C$.
Then $|\text{Aut}(C)| = 2$ when $g = 2$, 
and $|\text{Aut}(C)| = 1$ when $g > 2$. By (2.3), 
$$N_{g, d} = \pi^{-1}([C]) \cdot Z = \pi^{-1}([C_0]) \cdot Z
= {|\text{Aut}(C_0)| \over |\text{Aut}(C)|} 
\cdot \big (\overline{\frak M}_{C_0}(d) \cdot Z \big ).$$
Here $\overline{\frak M}_{C_0}(d)$ is identified with $\pi^{-1}([C_0])$
but with the reduced scheme structure. So to prove (1.1), 
it suffices to show that for $g = 2$ and $d \ge 4$, 
\smallskip
$$\overline{\frak M}_{C_0}(d) \cdot Z = {1 \over |\text{Aut}(C_0)|} \cdot 
\bigg[ {(d-1)(d-2)(d-3) \over d} N_d$$
$$+ \sum \limits_{d_1+d_2=d}
{d_1d_2(d_1d_2d-6d+18)-4d \over 6d}
{3d-2 \choose 3d_1 -1} d_1d_2 N_{d_1} N_{d_2} \bigg]. \tag 2.5$$

\smallskip\noindent
Note that $\overline{\frak M}_{C_0}(d) \cap Z \subset 
\overline{\frak M}_{C_0}(d) \cap \overline{W}_g(d)$. 
Let $[\mu: (D, p_1, \ldots, p_{3d - 2(g-1)})]$ be a point in 
$\overline{\frak M}_{C_0}(d) \cap \overline{W}_g(d)$.
Then $D$ consists of $(k+2(g-1))$ smooth rational curves with $k \ge 0$.
For simplicity, we also use $R_1, \ldots, R_{2(g-1)}$ to stand for 
the $2(g-1)$ smooth rational curves in $D$ which are identified with 
the $2(g-1)$ smooth rational curves $R_1, \ldots, R_{2(g-1)}$ in $C_0$
after $D$ is contracted to $C_0$.   
Dropping $R_1, \ldots, R_{2(g-1)}$ from $D$ results in a disjoint union 
$T_1 \coprod \ldots \coprod T_s$ of trees of smooth rational curves. 

\lemma{2.6} Assume $[\mu: (D, p_1, \ldots, p_{3d - 2(g-1)})]
\in \overline{\frak M}_{C_0}(d) \cap Z$.
Let $D_1, \ldots, D_m$ be all the irreducible components of $D$ 
such that $\mu|_{D_i}$ are not constant, 
and let $b_i = \text{deg}(\mu|_{D_i})$ for $1 \le i \le m$.
Then $m \le 2(g - 1)$; moreover, when $m = 2(g - 1)$,
$\mu(D_1), \ldots, \mu(D_m)$ have at most nodal singularities,
intersect each other transversally at nonsingular points, 
and have degrees $b_1, \ldots, b_m$ respectively.
\endproclaim
\noindent 
{\it Proof.} Note that $\sum_{i=1}^m b_i = d$.
For $1 \le i \le m$, let $\tilde b_i$ be the degree of $\mu(D_i)$ in $\Pee^2$. 
Then, $\tilde b_i \le b_i$. On the one hand, 
since $[\mu: (D, p_1, \ldots, p_{3d - 2(g-1)})] \in
\overline{\frak M}_{C_0}(d) \cap Z$, $\mu(D)$ has to pass 
$3d - 2(g-1)$ general points in $\Pee^2$. On the other hand,
the degree-$\tilde b_i$ irreducible rational curve $\mu(D_i)$ 
can pass through at most $(3 \tilde b_i - 1)$ general points in $\Pee^2$. 
So $\mu(D)$ can pass through at most
$\sum_{i=1}^m (3 \tilde b_i - 1) \le (3d -m)$ general points in $\Pee^2$.
Thus $m \le 2(g - 1)$. Moreover, if $m = 2(g - 1)$,
then $\tilde b_i = b_i$ and $\mu(D_1), \ldots, \mu(D_m)$ have 
at most nodal singularities and intersect each other transversally. \qed

\bigskip\noindent
{\it Proof of} (1.1): We shall now prove formula (1.1) by verifying (2.5).
So let $g = 2$ and $d \ge 4$. Put $d_i = \text{deg}(\mu|_{R_i})$ for 
$i = 1, 2$. Then, there are three cases:
\roster
\item"{(i)}" both $d_1$ and $d_2$ are positive;
\item"{(ii)}" exactly one of $d_1$ and $d_2$ is positive (so the other is zero);
\item"{(iii)}" $d_1 = d_2 = 0$.
\endroster
Our strategy is the following. Fix the nonnegative integer $k$. 
In each of the above three cases,  
we shall estimate the number $n(k)$ of moduli of various points 
$[\mu: (D, p_1, \ldots, p_{3d - 2})]$
in $\overline{\frak M}_{C_0}(d) \cap \overline{W}_2(d)$. 
Since the linear systems $e_i^*|\ell|$ are base-point-free,
it follows that if $n(k) < 6d - 4$, then the case will not contribute 
to the intersection number $\overline{\frak M}_{C_0}(d) \cap Z$.
We shall show that only cases (i) and (ii) with $k= 0$ may contribute. 
Furthermore, all cases (i) and (ii) with $k = 0$ actually contribute, 
i.e. any such map is actually in 
$\overline{\frak M}_{C_0}(d) \cap \overline{W}_2(d)$.
So formula (2.5) will be the sum of the contributions of 
cases (i) and (ii) with $k = 0$.  
Notice that by Lemma 2.6, we may assume that $m \le 2$ and 
that if $m = 2$, then $\mu(D_1)$ and $\mu(D_2)$ have 
at most nodal singularities and intersect transversally.

First of all, we consider case (i), that is, both $d_1$ and $d_2$ are positive.
By Lemma~2.6, $\mu$ is constant on the trees $T_i$.  
If $k=0$, then the number of moduli of these points 
$[\mu: (D, p_1, \ldots, p_{3d - 2})]$ in $\overline{\frak M}_{C_0}(d)$ is 
$$n(0) \le (3d - 2) + (3d - 2) = 6d - 4$$
where the first $(3d-2)$ is the number of moduli of the stable map $\mu$, 
and the second $(3d-2)$ is the number of moduli of the $(3d-2)$ marked points.
We claim that for arbitrary $k$, the number of moduli of these types of
maps satisfies $n(k)\le 6d-4-k$. To see this, consider the effect of
adding an additional rational component $D'$ to $D$ on which $\mu$ is constant.
If $D'$ meets the other components
in one point, it must contain two marked points for stability, and these
points do not have moduli.  Since $D'$ can replace at most one point that
has no moduli (the point $D'\cap D$), we have $n(k+1)\le n(k)-1$ for such
types of maps. The other possibility is for $D'$ to meet the other components
in two points. Then $D'$ must contain at least one point without
moduli, but no marked points have been replaced, so we obtain
$n(k+1)\le n(k)-1$ in this case as well. This proves the claim.
Thus only the cases with $k=0$ can occur, 
as claimed in the preceding paragraph. 
Furthermore, all cases (i) with $k = 0$ actually contribute.
Indeed, suppose that we have a map $\mu: C_0 (= D) \to \Pee^2$ 
with $d_1>0$ and $d_2>0$. Consider a general flat family of 
curves $\eta:\Cal E\to\Delta_t$ with $\eta^{-1}(0)=C_0$ and
$\eta^{-1}(t)$ smooth for $t\ne 0$. Consider the moduli functor of
relative line bundles on families of curves. The obstruction space for this
functor is 0, since it lies in an appropriate $\text{Ext}^2$ on a curve. 
So the functor is smooth. 
Since $C_0$ is stable, we can find a line bundle ${\Cal L}$ 
on $\Cal E$ which restricts to $\mu^*{\Cal O}_{\Pee^2}(1)$ on $C_0$ 
(possibly after a finite base change). Since $d>2$, 
either $d_1 \ge 2$ or $d_2 \ge 2$. 
By (2.4), $H^1(C_0, \mu^*{\Cal O}_{\Pee^2}(1)) = 0$. 
It follows that $R^1\eta_*{\Cal L}=0$ and that $\eta_*{\Cal L}$ is 
locally free of rank $(d-1)$. Thus there is no obstruction
to extending the three sections of $\mu^*{\Cal O}_{\Pee^2}(1)$ which
determine the map $\mu$ in the standard coordinates of $\Bbb P^2$.
This shows that $\mu$ is in $\overline{W}_2(d)$. So 
$\mu \in \overline{\frak M}_{C_0}(d) \cap \overline{W}_2(d)$.
Since $H^1(C_0, \mu^*T_{\Pee^2}) = 0$ for all such stable maps $\mu$ 
and the linear series $e_i^*|\ell|$ are base-point-free,
the contribution of case (i) to $\overline{\frak M}_{C_0}(d) \cap Z$
consists of finitely many reduced points. 
Notice that for a fixed pair of integers $d_1$ and $d_2$ 
with $d_1 > 0$ and $d_2 = d - d_1 > 0$, the number of unordered pairs 
$\{C_1, C_2\}$ of irreducible, reduced, nodal, degree-$d_1$ and degree-$d_2$ 
rational plane curves $C_1$ and $C_2$ whose union $C_1 \cup C_2$ 
pass through $(3d - 2)$ general points in $\Pee^2$ is 
${3d - 2 \choose 3d_1 -1} N_{d_1} N_{d_2}$. Thus the contribution of case (i)
to $\overline{\frak M}_{C_0}(d) \cdot Z$ is 
$${1 \over |\text{Aut}(C_0)|} \cdot \sum \limits_{d_1+d_2=d}
{d_1 d_2 \choose 3} {3d - 2 \choose 3d_1 -1} N_{d_1} N_{d_2}. \eqno (2.7)$$

For case (ii), we assume without loss of generality that 
$d_1 > 0$ and $d_2 = 0$. We start with $m = 2$. Let $\Pee$ be 
the unique smooth rational component of $D$ (aside from $R_1$) 
with nonconstant $\mu|_{\Pee}$. Let $T$ be the tree containing $\Pee$.
Then $T \cap R_1$ either is empty or consists of precisely one point.
If $T \cap R_1$ is empty, or if $T \cap R_1$ is nonempty but
the unique point in $T \cap R_1$ is not identified with one of 
the three singular points of $C_0$, 
then $\mu(R_1)$ has at least a triple point. 
If the unique point in $T \cap R_1$ is one of the singular points, 
then the other two singular points of $C_0$, 
when identified with points of $R_1$, 
are mapped by $\mu$ to a double point of $\mu(R_1)$, 
through which $\mu(\Pee)$ must pass. By Lemma 2.6,
the points $[\mu: (D, p_1, \ldots, p_{3d - 2})]$ can not be contained in 
$\overline{\frak M}_{C_0}(d) \cap Z$. So the case $m=2$ does not contribute.

Next, let $m = 1$. Then $\mu$ is constant on 
the closure $\overline{D \backslash R_1}$ of ${D \backslash R_1}$ in $D$.
So there exist at least three points $q_1, q_2, q_3$ in $R_1$
such that the images $\mu(q_1), \mu(q_2), \mu(q_3)$ are the same, 
i.e. every divisor in $(\mu|_{R_1})^*|\ell|$ 
which contains $p_1$ must contain both $p_2$ and $p_3$. 
This imposes $4$ independent conditions in choosing the linear series 
$(\mu|_{R_1})^*|\ell|$ from the complete linear system $|(\mu|_{R_1})^*\ell|$.
If $k=0$, the number of moduli of these points 
$[\mu: (D, p_1, \ldots, p_{3d - 2})]$
in $\overline{\frak M}_{C_0}(d) \cap \overline{W}_2(d)$ is at most
$$n(0) = [(3d -1) + (3-4)] + (3d-2) = 6d - 4.$$
Here, $[(3d -1) + (3-4)]$ is an upper bound for the number of moduli of
$\mu(R_1)$, where the integer $3$ in $(3-4)$ is the number of moduli of 
the three points $q_1, q_2, q_3$ varying in $R_1$. 
As in case (i) above, we see that for
general $k$, the number of moduli satisfies $n(k)\le 6d-4-k$.
It follows that only the case $k=0$ can occur.
Again by arguments  similar to those in case (i), 
we see that all cases (ii) with $k = 0$ actually contribute and 
that the contribution to $\overline{\frak M}_{C_0}(d) \cap Z$
consists of finitely many reduced points.
Furthermore, if $[\mu: (C_0, p_1, \ldots, p_{3d - 2})]$ is 
such a reduced point in $\overline{\frak M}_{C_0}(d) \cap Z$,
then $\mu(C_0) = \mu(R_1)$ is an irreducible and degree-$d$ rational 
plane curve that passes through $(3d - 2)$ general points in $\Pee^2$ and 
has at least one triple point. In fact, such a curve in $\Pee^2$
must have exactly one triple point with all other singularities being nodes.
Now there are only finitely many 
irreducible, reduced, degree-$d$ rational plane curves 
that pass through $(3d - 2)$ general points in $\Pee^2$ and 
have exactly one triple point with all other singularities being nodes.
The number $\tilde N_d$ of such curves is given by Lemma 3.2
which will be proved in the next section.
Taking into account of the automorphism group of $C_0$ and 
the symmetry between $R_1$ and $R_2$,
we see that the contribution of case (ii) to 
$\overline{\frak M}_{C_0}(d) \cdot Z$ is
$${1 \over |\text{Aut}(C_0)|} \cdot 2 \tilde N_d. \eqno (2.8)$$

Next we consider case (iii), that is, both $\mu|_{R_1}$ and $\mu|_{R_2}$
are constant. We start with $m = 1$. Let $\Pee$ be 
the unique smooth rational component of $D$ with nonconstant $\mu|_{\Pee}$,
and let $T$ be the unique tree in $T_1, \ldots, T_s$ 
such that $\Pee \subset T$. Here and in other subcases of case (iii), we
will be able  to assume without loss of generality that the tree $T$ is
actually a particularly simple chain. The recurring theme will be that if
$W$ is a subvariety of $\overline{\frak M}_{C_0}(d)$ 
such that $T$ is a certain type of chain 
and $\dim(W\cap\overline{W}_g(d))<6d-4$, then the same conclusion
will hold for subvarieties $W'$ associated to more complicated trees, since
in these situations $W'$ will always be a subvariety of the closure
$\overline{W}$ obtained from contracting suitable $\mu$-constant curves.
Returning to the subcase at hand, we may for the above reason assume that
$\Pee=T$ has one component.
If $T$ intersects $R_1$ or $R_2$ but not both,
then by Lemma 3.7 and Lemma 3.13, 
there exist $4$ independent conditions in choosing the linear series 
$(\mu|_T)^*|\ell|$ from $|(\mu|_T)^*\ell|$.
So the number of moduli of the points $[\mu: (D, p_1, \ldots, p_{3d - 2})]$
in $\overline{\frak M}_{C_0}(d) \cap \overline{W}_2(d)$ is at most
$$n(k) \le [(3d - k) + (2-4)] + (3d-2) = 6d - 4 - k \le 6d - 5$$
where the integer $2$ in $(2-4)$ is the number of the moduli of the point
$T \cap (R_1 \cup R_2)$ varying in both $T$ and $R_1 \cup R_2$. 
Thus this case makes no contribution to $\overline{\frak M}_{C_0}(d) \cap Z$.
If $T$ intersects both $R_1$ and $R_2$, then by Remark 3.12 (ii) and Lemma 3.7,
there still exist $4$ independent conditions in choosing  
$(\mu|_T)^*|\ell|$ from $|(\mu|_T)^*\ell|$.
So the number of moduli of the points $[\mu: (D, p_1, \ldots, p_{3d - 2})]$
in $\overline{\frak M}_{C_0}(d) \cap \overline{W}_2(d)$ is at most
$$n(k) \le [(3d - k) + (2-4)] + (3d-2) = 6d - 4 - k \le 6d - 5$$
where the integer $2$ in $(2-4)$ is the number of the moduli of the two points
$T \cap R_1$ and $T \cap R_2$ varying in $T$. It follows that 
this case makes no contribution to $\overline{\frak M}_{C_0}(d) \cap Z$.

We are left with case (iii) and $m = 2$. Then $k \ge 2$.
Let $D_1, D_2$ be the two rational components in $D$ such that 
$\mu|_{D_1}, \mu|_{D_2}$ are nonconstant. Then $D_i \ne R_j$ for $i, j = 1, 2$. 
First, assume that $D_1$ and $D_2$ are contained in the same tree $T$
from $T_1, \ldots, T_s$. 
If $T$ intersects $R_1$ or $R_2$ but not both, 
then by Lemma 3.7, there exist $4$ independent conditions in choosing  
the linear series $(\mu|_{D_1})^*|\ell|$ and $(\mu|_{D_2})^*|\ell|$ 
from the complete linear systems $|(\mu|_{D_1})^*\ell|$
and $|(\mu|_{D_2})^*\ell|$.
So the number of moduli of the points $[\mu: (D, p_1, \ldots, p_{3d - 2})]$
in $\overline{\frak M}_{C_0}(d) \cap \overline{W}_2(d)$ is at most
$$n(k) \le [(3d - k) + (2-4)] + (3d-2) = 6d - 4 - k \le 6d - 6$$
where the integer $2$ in $(2-4)$ is the number of the moduli of the point
$T \cap (R_1 \cup R_2)$ varying in both $T$ and $R_1 \cup R_2$,  
and this case makes no contribution to $\overline{\frak M}_{C_0}(d) \cap Z$.
Similarly, by Remark 3.12 (ii) and Lemma 3.7, 
the case when $T$ intersects both $R_1$ and $R_2$ makes no contribution 
to $\overline{\frak M}_{C_0}(d) \cap Z$.
Next, assume that $D_1$ and $D_2$ are contained in two different trees,
say $T_1$ and $T_2$, from $T_1, \ldots, T_s$.
If $T_1$ intersects both $R_1$ and $R_2$ and 
if $D_1$ intersects $\overline{D \backslash D_1}$ at least twice,
then $\mu(D_2)$ passes through a singular point of $\mu(D_1)$.
By Lemma 2.6, this case makes no contribution to 
$\overline{\frak M}_{C_0}(d) \cap Z$. So we may assume that 
if $T_i$ intersects both $R_1$ and $R_2$, 
then $D_i$ intersects $\overline{D \backslash D_i}$ once. 
By Lemma 3.7 and Remark 3.12 (ii), 
the number of moduli of the points $[\mu: (D, p_1, \ldots, p_{3d - 2})]$
in $\overline{\frak M}_{C_0}(d) \cap \overline{W}_2(d)$ is at most
$$n(k) \le [(3d - k) + (4-4)] + (3d-2) = 6d - 2 - k \le 6d - 4.$$
If $n(k) = 6d - 4$, then $k = s = 2$ and $T_i = D_i$ ($i = 1, 2$) 
intersects $R_1$ and $R_2$ but not both. Since the number of the moduli for 
the two points $D_1 \cap (R_1 \cup R_2)$ and $D_2 \cap (R_1 \cup R_2)$ 
varying in $R_1 \cup R_2$ is $2$, the pairs $(\mu(D_1), \mu(D_2))$
of curves form a codimension-$2$ subset in 
the $(3d-2)$-dimensional variety $S(0, d_1) \times S(0, d_2)$ where 
$S(0, d_i)$ stands for the moduli space of irreducible, reduced, nodal, 
degree-$d_i$ rational plane curves. It follows that this case makes 
no contribution to $\overline{\frak M}_{C_0}(d) \cap Z$. 

Summing up (2.7) and (2.8) and using (3.3), 
we obtain (2.5) and hence (1.1). \qed

\bigskip
Finally, some remarks about the number $N_{g, d}$ with $g > 2$ follow. 
We expect that only the cases with $k = 0$ (that is, stable maps
$[\mu: (D, p_1, \ldots, p_{3d - 2(g-1)})]$ with $D = C_0 = C_{0, g}$) 
contribute to  $N_{g, d} = |\text{Aut}(C_0)| \cdot 
\big (\overline{\frak M}_{C_0}(d) \cdot Z \big )$. 
For instance, as in the proof of the above case (i), 
the case when $k = 0$ and $\text{deg}(\mu|_{R_i}) > 0$ for 
all $1 \le i \le 2(g-1)$ contributes to $N_{g, d}$, and its contribution is 
$$(3d - 2(g-1))! \cdot \sum \limits_{d_1 + \ldots + d_{2(g-1)} = d} \,
\prod \limits_{i=1}^{2(g-1)} {d_i^3 \cdot N_{d_i} \over (3d_i -1)!}.
\eqno (2.9)$$
It is interesting to notice that for $g>2$, the term~(2.9) is precisely
a term that arises in the calculation of $\Psi_{g,d}$.
In general, let $S$ be any nonempty subset of $\{1, \ldots, 2(g-1) \}$.
We believe that the contribution of the cases with $k = 0$ and 
$\text{deg}(\mu|_{R_i}) > 0$ if and only if $i \in S$ is closely related
to the number of irreducible, reduced, degree-$d$ plane rational curves
which pass through certain general points in $\Pee^2$ and 
have certain types of singularities. 

\bigskip\noindent
{\bf 3. Rational plane curves with triple points and 
constraint on stable maps}

In this section, we prove the lemmas quoted in the previous section.
We shall compute the number $\tilde N_d$ of 
irreducible, reduced, degree-$d$ rational plane curves 
that pass through $(3d-2)$ general points in $\Pee^2$ and 
have exactly one triple point with all other singularities being nodes. 
Then we analyze the constraints on stable maps that are degenerations 
of maps whose domains are irreducible and smooth curves. 

\ssection{3.1. Rational plane curves with triple points}

Fix $d \ge 3$. Intersections of $\Bbb Q$-divisors in the moduli space
$\overline{\frak M}_{0, 0}(\Bbb P^2, d)$ has been studied by Pandharipande.
We recall some results from \cite{Pa2, Pa3}.
The boundary $\Delta$ of $\overline{\frak M}_{0, 0}(\Bbb P^2, d)$ is 
the locus corresponding to stable maps whose domains are reducible,
and is of pure codimension-$1$.
For $1 \le i \le \left [ {d \over 2} \right ]$, 
let $K^i$ be the irreducible component of $\Delta$ whose general elements
are the form $\mu: D_1 \cup D_2 \to \Pee^2$ such that 
$\text{deg}(\mu|_{D_1}) = i$, and $D_1$ and $D_2$ are smooth rational curves
meeting transversely at one point. 
Let $\Cal H$ be the locus of $\overline{\frak M}_{0, 0}(\Bbb P^2, d)$
corresponding to maps whose images pass through a fixed point in $\Pee^2$.
Then $\Cal H$ is a Cartier divisor in $\overline{\frak M}_{0, 0}(\Bbb P^2, d)$,
and $\text{Pic}(\overline{\frak M}_{0, 0}(\Bbb P^2, d)) \otimes \Bbb Q$
is generated by $\Cal H$ and $K^i$ with 
$1 \le i \le \left [ {d \over 2} \right ]$; moreover,
$$\Cal H^{3d-1} = N_d, \quad K^i \cdot \Cal H^{3d-2} = 
\cases 
{3d-2 \choose 3i -1} i(d-i) N_i N_{d-i}, &\text{if $i \ne {d \over 2}$}\\
{1 \over 2} {3d-2 \choose 3{d \over 2} -1} \left ({d \over 2} \right )^2 
N_{d \over 2}^2, &\text{if $i = {d \over 2}$}.
\endcases  \tag 3.1$$
Let $Z \subset \overline{\frak M}_{0, 0}(\Bbb P^2, d)$ 
be the subvariety consisting of degree-$d$ maps $\mu: \Pee^1 \to \Pee^2$
such that $\mu(\Pee^1)$ has exactly one triple point 
with all other singularities being nodes, and let ${\overline Z}$ be 
the closure of $Z$ in $\overline{\frak M}_{0, 0}(\Bbb P^2, d)$.
Then $Z$ is reduced and of pure codimension-$1$, 
and ${\overline Z}$ is a Weil divisor.
Similar arguments as in section (3.4) of \cite{Pa2} show that 
the intersection ${\overline Z} \cap \Cal H^{3d-2}$ determined by $(3d-2)$ 
general points in $\Pee^2$ consists of finitely many reduced points in $Z$. 
Thus $\tilde N_d = {\overline Z} \cdot \Cal H^{3d-2}$.

\lemma{3.2} For $d \ge 3$, $\tilde N_d$ can be expressed in terms of the $N_d$:
$${(d-1)(d-2)(d-3) \over 2d} N_d - 
\sum_{d_1+d_2=d} {d_1 d_2 (d-6)+2d \over 4d}
{3d-2 \choose 3d_1 -1} d_1 d_2 N_{d_1} N_{d_2}. \eqno (3.3)$$
\endproclaim
\noindent
{\it Proof.} 
Let ${\overline Z} = a \Cal H + \sum_{i = 1}^{\left [ {d \over 2} \right ]} 
a_i K^i$ where $a, a_i \in \Bbb Q$. 
By (3.1), it suffices to determine $a$ and $a_i$.
Fix a nonsingular curve $C$. As in section (4.5) of \cite{Pa2}, 
we compute the intersection number $C \cdot \lambda^*{\overline Z}$ 
for some morphism $\lambda: C \to \overline{\frak M}_{0, 0}(\Bbb P^2, d)$ 
which is constructed as follows.
Let $\pi_2: S = \Pee^1 \times C \to C$ be the second projection,
and let $\Cal N$ be a line bundle on $S$ of degree type $(d, k)$ 
where $k$ is a very large integer. Let $z_0, z_1, z_2 \in H^0(S, \Cal N)$
determine a rational map $\phi: S \dashrightarrow \Pee^2$ such that
\roster
\item"{(i)}" every base point of $\phi$ is simple. Here a base point $s$ 
of $\phi$ is {\it simple of degree $i$ with $1 \le i \le d$} if 
the blowing-up of $S$ at $s$ resolves $\phi$ locally at $s$ and 
the resulting map is of degree $i$ on the exceptional divisor;
\item"{(ii)}" there exist only finitely many points $c_1, \ldots, c_n \in C$
such that for each $i$, $c_i$ is not the projection of base points of $\phi$
and $\overline \phi(\overline \pi^{-1}(c_i))$ contains 
one and exactly one triple point. Here $\overline S$ is the blow-up of 
$S$ at the base points, $\overline \pi: \overline S \to C$ is 
the projection to $C$, and $\overline \phi: \overline S \to \Pee^2$ is 
the resolution of $\phi$.
\endroster 
Now $\overline \phi: \overline S \to \Pee^2$ and 
$\overline \pi: \overline S \to C$ induce a morphism 
$\lambda: C \to \overline{\frak M}_{0, 0}(\Bbb P^2, d)$.
A point of the intersection $C \cdot \lambda^*K^i$ can arise in two cases,
that is, a simple base point of degree $i$ or $(d-i)$ can be blown-up.
Let $C \cdot \lambda^*K^i = x_i + y_i$ where $x_i$ and $y_i$ are 
the number of instances of the first and second case respectively. 
On the one hand,  
$C \cdot \lambda^*\Cal H = 2dk - \sum_{i= 1}^{\left [ {d \over 2} \right ]} 
[i^2x_i + (d-i)^2y_i]$ according to \cite{Pa2}. Thus,
$$C \cdot \lambda^*{\overline Z} = 
2adk - \sum_{i= 1}^{\left [ {d \over 2} \right ]} 
a[i^2x_i + (d-i)^2y_i] + \sum_{i = 1}^{\left [ {d \over 2} \right ]} a_i
(x_i + y_i). \eqno (3.4)$$
On the other hand, the triple-point formula of \cite{Kle} can be applied to
the map $(\overline{\phi},\overline{\pi}):\overline{S}\to \Pee^2\times C$ to yield
$$C \cdot \lambda^*{\overline Z} = (d-1)(d-2)(d-3)k
+ \sum_{i= 1}^{\left [ {d \over 2} \right ]} (-{1\over 2} i^2 d^2
+3i^2d -{1\over 2} id-1+3i-5i^2) x_i$$ 
$$+\sum_{i= 1}^{\left [ {d \over 2} \right ]} (-{1\over 2} (d-i)^2 d^2
+3(d-i)^2d -{1\over 2} (d-i)d-1+3(d-i)-5(d-i)^2) y_i. \eqno (3.5)$$
Parts of this computation were performed using
Schubert \cite{K-S}.
Comparing the coefficients of $k, x_i, y_i$ in (3.4) and (3.5) leads to
$$a = {(d-1)(d-2)(d-3) \over 2d} \qquad \text{and} \qquad
a_i = -{i(d-i)(d-6)+2d \over 2d}.  \eqno (3.6)$$
Therefore the formula for $\tilde N_d$ follows from (3.1), (3.6) and 
$$\tilde N_d = {\overline Z} \cdot \Cal H^{3d-2} = a \Cal H^{3d-1} + 
\sum_{i = 1}^{\left [ {d \over 2} \right ]} a_i K^i \cdot \Cal H^{3d-2}.
\qed$$

\ssection{3.2. Constraint on stable maps}

Let $\eta: \Cal E \to \Delta_t$ be a flat family of stable genus-$g$ curves 
in $\overline{\frak M}_g$ satisfying 
\roster
\item"{(a)}" $\eta^{-1}(t)$ is irreducible and smooth for all $t \ne 0$.
\item"{(b)}" for every smooth point $p \in C_0 {\overset \text{def} \to =}
\eta^{-1}(0)$, there exists a basis $\Lambda_1, \ldots, \Lambda_g$
for $H^0(C_0, \omega_{C_0})$ such that locally near $p$, we have
$\Lambda_j = v^{j-1} f_j(v) \cdot
\text{d}v$ for some holomorphic functions $f_j(v)$ satisfying $f_j(0) \ne 0$ 
for $1 \le j \le g$ where $v$ is a local coordinate of $C_0$ centered at $p$.
\endroster
Let $\hat \eta: \hat \Cal E \to \Delta_t$ be the family obtained by 
blowing-up $\Cal E$ at points in $C_0$ (possibly infinitely near) 
and by adding $(3d-2(g-1))$ markings. Let $D = \hat \eta^{-1}(0)$,
and $T_1, \ldots, T_s$ be all the connected components of 
the closure $\overline{D \backslash C_0}$ of $D \backslash C_0$ in $D$.
Then each $T_i$ is a tree of smooth rational curves, and 
$D = C_0 \cup (\coprod_{i=1}^s T_i)$.
Assume that $\mu: \hat \Cal E \to \Pee^2$ is a morphism such that 
$\mu, \hat \eta$, and the $(3d-2(g-1))$ markings determine
a family of stable maps in $\overline{W}_g(d)$.
Furthermore, suppose that there exist smooth rational components 
$D_1, \ldots, D_m$ contained in the trees $T_1, \ldots, T_s$
with $\text{deg}(\mu|_{D_i}) > 0$ for each $i$ and 
$\sum_{i = 1}^m \text{deg}(\mu|_{D_i}) = d$.

\lemma{3.7} Let $g \ge 2$ and $d > 2(g-1)$. 
Assume that $\hat \Cal E$ is obtained from $\Cal E$ by a chain of blowups 
at smooth points of $C_0$ (possibly infinitely near). 
If $m \le 2$, then there exist $2g$ independent conditions
in choosing the $m$ linear series $(\mu|_{D_i})^*|\ell|$, $1 \le i \le m$
from the $m$ complete linear systems $|(\mu|_{D_i})^*\ell|$, $1 \le i \le m$.
\endproclaim
\proof
There are three separate cases: (i) $m = 1$; (ii) $m = 2$, and $D_1, D_2$
are contained in the same tree in $T_1, \ldots, T_s$; 
(iii) $m = 2$, and $D_1, D_2$
are contained in the two different trees in $T_1, \ldots, T_s$.
The proofs of (ii) and (iii) are very similar to the proof of (i) 
but need some extra preparation.

{\bf Case (i)}: We follow the approach in \cite{Pa1}. 
For simplicity, we first assume that 
$\hat \Cal E$ is the blow-up of $\Cal E$ at a smooth point $p \in C_0$.
Then $s=1= i$, and $\Pee = T_1 = \Pee^1$ is the exceptional divisor. 
For each $j$ with $1 \le j \le d$, let $\Cal G_j = \Cal H_j \subset \Cal E$
be the open subset of $\Cal E$ on which the morphism $\eta$ is smooth. Put
$$X = \Cal G_1 \times_{\Delta_t} \ldots \times_{\Delta_t} \Cal G_d 
\times \Cal H_1 \times_{\Delta_t} \ldots \times_{\Delta_t} \Cal H_d.$$
Then $X$ is a smooth open subset of the $2d$-fold fiber product of $\Cal E$
over $\Delta_t$. Let $Y \subset X$ be the subset of points
$y = (g_1, \ldots, g_d, h_1, \ldots, h_d)$ such that the two divisors
$\sum_{j} g_j$ and $\sum_{j} h_j$ are linearly equivalent on the curve
$\eta^{-1}(\eta(y))$. 

Let $\gamma: \Delta_t \to \Cal E$ be a local holomorphic section of $\eta$ 
such that $\gamma(0) = p$. Let $V$ be a nonvanishing local holomorphic field of
vertical tangent vectors to $\Cal E$ on an open subset containing $p$. 
The section $\gamma$ and the vertical vector field $V$ together determine 
local holomorphic coordinates $(t, v)$ on $\Cal E$ at $p$. 
Let $\phi_V: \Cal E \times \Cee \to \Cal E$ be the holomorphic flow of $V$
defined locally near $(p, 0) \in \Cal E \times \Cee$. Then the coordinate map
$\psi: \Cee^2 \to \Cal E$ is defined by $\psi(t, v) = \phi_V(\gamma(t),v)$.
Since $p \in C_0$ is a smooth point in $C_0$, $p \in \Cal G_j$ and 
$p \in \Cal H_j$ for $1 \le j \le d$. Put $x_p = 
(p, \ldots, p, p, \ldots, p) \in X$. Then the local coordinates on $X$ 
near $x_p$ are given by $(t, v_1, \ldots, v_d, w_1, \ldots, w_d)$, 
and the coordinate map $\psi_X$ is defined by putting 
$\psi_X(t, v_1, \ldots, v_d, w_1, \ldots, w_d)$ to be 
$$(\psi(t, v_1), \ldots, \psi(t, v_d), \psi(t, w_1), \ldots, \psi(t, w_d))
\in X.$$
By our assumption, there exists a basis $\Lambda_1, \ldots, \Lambda_g$
of $H^0(C_0, \omega_{C_0})$ such that near $p$, we have 
$\Lambda_k(v) = v^{k-1} f_k(v) \cdot
\text{d}v$ and $f_k(0) \ne 0$ for $1 \le k \le g$.  
Let $\Lambda_k(t, v)$ ($1 \le k \le g$) be a local holomorphic extension 
of $\Lambda_k(v)$ at $0 \in \Delta_t$ such that for a fixed $t$, 
$\Lambda_1(t, v), \ldots, \Lambda_g(t, v)$
form a basis of $H^0(\eta^{-1}(t), \omega_{\eta^{-1}(t)})$ using the
coordinate map $\psi$.  This basis may be chosen so that 
$$\Lambda_k(t, v) = v^{k-1} f_k(t, v) \cdot \text{d}v$$
with $f_k(0, 0) = f_k(0) \ne 0$. Now the local equations of $Y \subset X$ 
at the point $x_p$ are
$$\sum_{j=1}^d \left ( \int_0^{v_j} \Lambda_k(t, v) - 
\int_0^{w_j} \Lambda_k(t, v) \right ) = 0, \qquad 1 \le k \le g.
\eqno (3.8)$$

Let $L_1$ and $L_2$ be general divisors in $\mu^*|\ell|$ such that 
each intersects $\Pee$ at $d$ distinct points. 
For $1 \le \alpha \le 2$, $L_\alpha$ determines local holomorphic sections
$s_{\alpha, 1} + \ldots + s_{\alpha, d}$ of $\hat \eta$ at $0 \in \Delta_t$.
These sections $s_{\alpha, j}$ with $1 \le \alpha \le 2$ and $1 \le j \le d$
determine a map $\lambda: \Delta_t \to Y$ locally at $0 \in \Delta_t$.
Let an affine coordinate on $\Pee^1$ be given by $\xi$ corresponding to 
the normal direction 
$${\text{d}\gamma \over \text{d}t}|_{t = 0} + \xi \cdot V(p).$$
Let $s_{1, j}(0) = \nu_j \in \Cee^1 \subset \Pee^1$ and 
$s_{2, j}(0) = \omega_j \in \Cee^1 \subset \Pee^1$ be given in terms of 
the affine coordinates $\xi$. Then the map $\lambda$ has the form
$$\lambda(t) = (t, \nu_1(t), \ldots, \nu_d(t), 
\omega_1(t), \ldots, \omega_d(t)) \eqno (3.9)$$ 
where 
$$\nu_j(t) = \nu_j t + O(t^2), \ \omega_j(t) = \omega_j t + O(t^2),\qquad 
1 \le j \le d. \eqno (3.10)$$

Since $Y$ is defined by the equations (3.8), we obtain 
$$\sum_{j=1}^d \left ( \int_0^{\nu_j(t)} v^{k-1} f_k(t, v) \cdot \text{d}v - 
\int_0^{\omega_j(t)} v^{k-1} f_k(t, v) \cdot \text{d}v \right ) = 0, 
\quad 1 \le k \le g.  \eqno (3.11)$$
Differentiating (3.11) $k$-times with respect to $t$
and evaluating at $t = 0$ results in
$$\sum_{j=1}^d \left ( (k-1)! f_k(0,0) \cdot \nu_j^k - 
(k-1)! f_k(0,0) \cdot \omega_j^k \right ) = 0, \qquad 1 \le k \le g.$$
Since $f_k(0, 0) = f_k(0) \ne 0$, we have 
$\sum_{j=1}^d \nu_j^k = \sum_{j=1}^d \omega_j^k$ where $1 \le k \le g$.
Let $\beta_k$ be the $k$-th elementary symmetric function in $d$ variables.
Then, $\beta_k(\nu_1, \ldots, \nu_d) = 
\beta_k(\omega_1, \ldots, \omega_d)$ for $1 \le k \le g$. 
Put $\beta_k' = (-1)^k \cdot \beta_k(\nu_1, \ldots, \nu_d)$
for $1 \le k \le g$. Then the divisors in $(\mu|_{\Pee})^*|\ell|$ 
correspond to degree-$d$ polynomials of the form
$$K(\xi^d + \beta_1' \cdot \xi^{d-1} + \ldots + 
\beta_g' \cdot \xi^{d-g} + \ldots)$$
where $K$ stands for constants. It follows that the linear series
$(\mu|_{\Pee})^*|\ell|$ has vanishing sequence $\{0, \ge (g+1), *\}$ at
the point $\xi = \infty$ which is the intersection $C_0 \cap \Bbb P$.
Since the complete linear system $|(\mu|_{\Pee})^*\ell|$ is base-point-free,
the existence of a vanishing sequence of the form $\{0, \ge (g+1), *\}$
for the linear series $(\mu|_{\Pee})^*|\ell|$
imposes $2g$ independent conditions in choosing  
$(\mu|_{\Pee})^*|\ell|$ from $|(\mu|_{\Pee})^*\ell|$.

The general case arises when $n$ blowups are needed to obtain $\Pee$.  
In this situation, using automorphisms of the rational components in $T_i$,
we may assume that the form of $\lambda$ is again given by~(3.9) with~(3.10) 
replaced by
$$\nu_j(t) = \nu_j t^n + O(t^{n+1}), \,
\omega_j(t) = \omega_j t^n + O(t^{n+1}), \qquad 1 \le j \le d.$$
Then the calculation concludes as before. (Compare with \cite{Pa1}.)

{\bf Case (ii)}: For simplicity, we assume that $\hat \Cal E$ is 
the $2$-fold blow-up of $\Cal E$ at a smooth point $p \in C_0$.
Then $s=1= i$, and $T_1 = D_1 \cup D_2$ is the union of 
the two exceptional divisors. 
Let $d_i = \text{deg}(\mu|_{D_i})$ for $i = 1, 2$. 
Then, $d_1 > 0$, $d_2 > 0$, and $d_1 + d_2 = d$.
Let $L_1$ and $L_2$ be general divisors in $\mu^*|\ell|$ such that 
each intersects $D_i$ at $d_i$ distinct points. 
Let other notations be as in Case (i). 
Then locally at $0 \in \Delta_t$, $L_1$ and $L_2$ induce a map 
$\lambda: \Delta_t \to Y$ sending $t \in \Delta_t$ to 
the following point in $Y$:
$$\lambda(t) = (t, \nu_{1,1}(t), \ldots, \nu_{1, d_1}(t), 
\nu_{2,1}(t), \ldots, \nu_{2, d_2}(t),$$
$$\quad \omega_{1,1}(t), \ldots, \omega_{1, d_1}(t),
\omega_{2,1}(t), \ldots, \omega_{2, d_1}(t))$$ 
where for $1 \le i \le 2$ and $1 \le j \le d_i$, 
we may assume that $\nu_{i, j}(t) = \nu_{i, j} t^i + O(t^{i+1})$ and  
$\omega_{i, j}(t) = \omega_{i, j} t^i + O(t^{i+1})$ 
for some constants $\nu_{i, j}$ and $\omega_{i, j}$.
A similar argument as in Case (i) shows that
the linear series $(\mu|_{D_1})^*|\ell|$ has vanishing sequence 
$\{0, \ge (g+1), *\}$ at the point $C_0 \cap D_1$.
So there exist $2g$ independent conditions
in choosing $(\mu|_{D_1})^*|\ell|$ from the complete linear system 
$|(\mu|_{D_1})^*\ell|$.

{\bf Case (iii)}: Again for simplicity, we assume that $\hat \Cal E$ is 
the $2$-fold blow-up of $\Cal E$ at two smooth points $p_1, p_2 \in C_0$.
Then $s=2$, and $\{ T_1, T_2 \} = \{D_1, D_2\}$ is the set of  
the two exceptional divisors. 
Let $d_i = \text{deg}(\mu|_{D_i})$ for $i = 1, 2$. 
Then, $d_1 > 0$, $d_2 > 0$, and $d_1 + d_2 = d$.
Since $d > 2(g-1)$, we may assume that $d_1 \ge g$.
As in Case (i), we construct local coordinates $(t, v_i)$ and 
coordinate map $\psi_i$ on $\Cal E$ at each point $p_i$.
Let $X$ and $Y$ be as in Case (i). Define a point $x_{p_1,p_2} \in X$ by 
$$x_{p_1,p_2} = (\underbrace{p_1, \ldots, p_1}_{\text {$d_1$ times}}, 
\underbrace{p_2, \ldots, p_2}_{\text {$d_2$ times}},
\underbrace{p_1, \ldots, p_1}_{\text {$d_1$ times}}, 
\underbrace{p_2, \ldots, p_2}_{\text {$d_2$ times}}).$$
The local coordinates on $X$ near $x_{p_1,p_2}$ are given by 
$(t, v_1, \ldots, v_d, w_1, \ldots, w_d)$, 
and the coordinate map $\psi_X$ is defined by putting 
$\psi_X(t, v_1, \ldots, v_d, w_1, \ldots, w_d)$ to be 
$$(\psi_1(t, v_1), \ldots, \psi_1(t, v_{d_1}), 
\psi_2(t, v_{d_1 +1}), \ldots, \psi_2(t, v_d),$$
$$\psi_1(t, w_1), \ldots, \psi_1(t, w_{d_1}), 
\psi_2(t, w_{d_1 +1}), \ldots, \psi_2(t, w_d)).$$
Note that $x_{p_1,p_2} \in Y$ and the local equations of $Y \subset X$
at $x_{p_1,p_2}$ are given by (3.8) where 
$\Lambda_1(t, v), \ldots, \Lambda_g(t, v)$ are chosen so that 
$\Lambda_k(t, v_1) = v_1^{k-1} f_k(t, v_1) \cdot \text{d}v_1$
with $f_k(0, 0) \ne 0$. Let $L_1$ and $L_2$ be general divisors 
in $\mu^*|\ell|$ such that each intersects $D_i$ at $d_i$ distinct points.
Then locally at $0 \in \Delta_t$, $L_1$ and $L_2$ induce  
$\lambda: \Delta_t \to Y$ by
$$\lambda(t) = (t, \nu_{1,1}(t), \ldots, \nu_{1, d_1}(t), 
\nu_{2,1}(t), \ldots, \nu_{2, d_2}(t),$$
$$\quad \omega_{1,1}(t), \ldots, \omega_{1, d_1}(t),
\omega_{2,1}(t), \ldots, \omega_{2, d_1}(t))$$ 
where for $1 \le i \le 2$ and $1 \le j \le d_i$, 
we have $\nu_{i, j}(t) = \nu_{i, j} t + O(t^2)$ and  
$\omega_{i, j}(t) = \omega_{i, j} t + O(t^2)$ 
for some constants $\nu_{i, j}$ and $\omega_{i, j}$.
Now a similar argument as in Case (i) shows that for $1 \le k \le g$, 
$\sum_{j=1}^{d_1} \nu_{1,j}^k - \sum_{j=1}^{d_1} \omega_{1, j}^k$ is a homogeneous polynomial in $\nu_{2,1}, \ldots, \nu_{2,d_2}$
and in $\omega_{2,1}, \ldots, \omega_{2,d_2}$. 
Therefore for a fixed linear series $(\mu|_{D_2})^*|\ell|$ 
in $|(\mu|_{D_2})^*\ell|$, these exist $2g$ independent conditions
in choosing $(\mu|_{D_1})^*|\ell|$ from $|(\mu|_{D_1})^*\ell|$,
i.e. there exist $2g$ independent conditions
in choosing the linear series $(\mu|_{D_1})^*|\ell|$ and 
$(\mu|_{D_2})^*|\ell|$ from the complete linear systems 
$|(\mu|_{D_1})^*\ell|$ and $|(\mu|_{D_2})^*\ell|$.

This completes the proof of Case (iii) and hence the proof of the lemma.
\endproof

\noindent
{\it Remark 3.12.} (i) It is reasonable to expect that Lemma 3.7 holds 
for any $m$.
\par
(ii) In Lemma 3.7, we have assumed that $\hat \Cal E$ is the blow-up
of $\Cal E$ at smooth points of $C_0$. However, a slight modification of 
its proof shows that the conclusion is still true if $g = 2$ and 
the blowing-up $\hat \Cal E \to \Cal E$ also takes place 
at nodal points of $C_0$.

Finally, we show that the stable genus-$2$ curve $C_0$ constructed in 
Section~2 satisfies hypothesis (b) leading up to the statement of Lemma~3.7.

\lemma{3.13} For every smooth point $p \in C_0$, there exists a basis 
$\Lambda_1, \Lambda_2$ for $H^0(C_0, \omega_{C_0})$ such that 
$\Lambda_j = v^{j-1} f_j(v) \cdot \text{d}v$ and some holomorphic
functions $f_j(v)$ such that $f_j(0) \ne 0$ for 
$1 \le j \le 2$ where $v$ is a local coordinate of $C_0$ centered at $p$.
\endproclaim
\noindent
{\it Proof.}
Assume that $p \in R_1 = \Pee^1$. Choose an affine coordinate $z$ for $R_1$
such that the three nodal points in $C_0$ are identified with 
$0, 1, \infty$ in $R_1$. Then a basis for $H^0(C_0, \omega_{C_0})$
can be identified with $\Lambda_1' = {1 \over z} \cdot \text{d}z, 
\Lambda_2' = {1 \over z-1} \cdot \text{d}z.$
Let $z_0$ be the coordinate of $p \in R_1$, and let $v = z-z_0$. 
Then the desired basis consists of
$$\Lambda_1 = {1 \over z} \cdot \text{d}z 
= {1 \over v+z_0} \cdot \text{d}v, \quad
\Lambda_2 = {(z-z_0) \over z(z-1)} \cdot \text{d}z = 
{v \over (v+z_0)(v+z_0-1)} \cdot \text{d}v.  \qed$$

\enddocument